\documentclass[10pt,twocolumn,twoside]{IEEEtran}

\usepackage{graphicx} 
\usepackage[usenames,dvipsnames,table]{xcolor}
\usepackage[colorlinks=true, allcolors=black, urlcolor=darkgray]{hyperref}
\usepackage{url}

\usepackage{color}

\usepackage{amsmath, amssymb}  

\newcommand{\ones} {\mathbf{1}}
\newcommand{\R} {\mathbb{R}}
\newcommand{\T} {^\top}
\newcommand{\diag} {\mathrm{diag}}

\newcommand{\gram}[1]{\mathcal{W}_{\mathcal{K},#1}}
\newcommand{\Co}[1]{\mathcal{C}_{#1}}
\newcommand{\Cob}[1]{\bar{\mathcal{C}}_{#1}}
\newcommand{\lmin}[1]{\lambda_{\min}(\mathcal{W}_{\mathcal{K},#1})}

\newcommand{\lminhat}[1]{\hat{\lambda}_{\min}(\mathcal{W}_{\mathcal{K},#1})}
\newcommand{\inv} {^{-1}}
\newcommand{\N} {\mathbb{N}}
\newcommand{\G}{\mathcal{G}}

\newcommand{\E}{\mathcal{E}}
\newcommand{\V}{\mathcal{V}}
\newcommand{\K}{\mathcal{K}}
\newcommand{\eps}{\varepsilon}

\newtheorem{lem}{Lemma}
\newtheorem{thm}{Theorem}
\newtheorem{rem}{Remark}
\newtheorem{cor}{Corollary}
\newtheorem{prop}{Proposition}
\newtheorem{defn}{Definition}

\newtheorem{exmp}{Example}
\newtheorem{assum}{Assumption}

\def\Rw{\color{black}}

\begin{document}

\title{On the Role of Network Centrality in the Controllability of Complex Networks}

\author{Nicoletta~Bof,
        Giacomo~Baggio, 
        and~Sandro~Zampieri
\thanks{N.~Bof, G.~Baggio and S.~Zampieri are with the Dipartimento di Ingegneria dell'Informazione,
         Universit\`a di Padova,
         via Gradenigo, 6/B -- I-35131 Padova, Italy.
        E-mails: {\tt \{bofnicol, baggio, zampi\}@dei.unipd.it}.}}

\maketitle

\begin{abstract}
In recent years complex networks have gained increasing attention in different fields of science and engineering. The problem of controlling these networks is an interesting and challenging problem to investigate. In this paper we look at the controllability problem focusing on the energy needed for the control. Precisely not only we want to analyze whether a network can be controlled, but we also want to establish whether the control can be performed using a limited amount of energy. We restrict our study to irreducible and (marginally) stable networks and we find that the leading right and left eigenvectors of the network matrix play a crucial role in this analysis. Interestingly, our results suggest the existence of a connection between controllability and network centrality, a well-known concept in network science. In case the network is reversible, the latter connection involves the PageRank, an extensively studied type of centrality measure. Finally, the proposed results are applied to examples concerning random graphs.
\end{abstract}


%

\section{Introduction}
Complex networks are systems composed of a large number of units which interact
among themselves, forming in this way a behavior which is much richer than the behavior of the
single units \cite{Strogatz:01}.
Many systems, which model both natural processes and engineering structures, can be seen as complex networks. Among them, one can mention genomic networks and ecologic networks in biology, social networks  in sociology and economic or financial networks in economics, while in engineering, electric power grids, transportation networks and communication networks are some important examples \cite{Newman:10}.
Therefore many areas of science and technology can take a great advantage from a deep understanding of this class of systems.

Consequently, the properties of complex networks have attracted a lot of interest among different scientific communities in the last years. Controllability is one of these properties, and it consists in the possibility of steering the state of the network from any initial value into a final arbitrary one, by fixing the profile of the state of a subset of nodes, called control nodes, which are assumed to be directly accessible by the controller \cite{Kailath:80}. 
This property has been analyzed by relating it to the features of the underlying graph \cite{Liu:11}. In particular, the results on structural controllability theory \cite{Lin:74} have been exploited. This kind of approach  aims to establish whether the network is controllable and, specifically, which nodes must be controlled in order to obtain the network controllability.

From \cite{Liu:11}, several other papers focusing on this type of problem have been proposed, including \cite{Olshevsky:14,Olshevsky:15, Pequito:14, Cowan:12}.
However, an important aspect is neglected in this line of research. In fact, one can notice that, even if
a control profile able to drive the network state may exist in principle, this profile may not be physically implementable due to the energy it requires.
This observation gave rise to another type of approach related to network controllability, focusing specifically on the evaluation of a ``degree'' of controllability of the network.
This concept can be made precise only by introducing a suitable metric related to the energy needed for the control.
The evaluation of the control energy naturally involves, for linear systems, the notion of controllability Gramian,
which is a symmetric positive definite matrix in case the network is controllable \cite{Kailath:80}. 
The analysis of the controllability Gramian offers different ways to study the degree of controllability of a network, depending on which property of this matrix is inspected.
For instance, one can select the minimum eigenvalue of the controllability Gramian \cite{Yan:12}, \cite{PZB:14}, the trace of its inverse \cite{Summers:14}, its determinant \cite{Muller:72}, or its condition number \cite{Sun:13}.

In this context, different types of problems can be addressed. For example, one interesting issue consists in the optimal placement of the control nodes in a way such that the degree of controllability is maximized \cite{Summers:14, Summers:14-CDC,Tzoumas:14, Tzoumas:15-ACC}. Another engaging problem is to relate the network structure to the number of control nodes needed to make the network practically controllable, namely controllable by control profiles with bounded energy \cite{Yan:12, Sun:13, PZB:14, PZ-CDC:14, Gang-Yan:15, Enyioha:14}.

In this paper, we focus on this second problem. Precisely, we study in which cases the energy needed for the control tends to infinity as the number of nodes gets larger and larger. This can be a way to establish the classes of complex networks which are practically impossible to control, since the energy they require for their control can go over any threshold as the number of nodes increases. This behavior will depend both on the properties of the network and on the number of control nodes. As in \cite{PZB:14} and \cite{PZ-CDC:14}, we adopt the minimum eigenvalue of the controllability Gramian as the measure of the controllability degree. In our analysis we will restrict to networks with non-negative weight matrices which are irreducible and (marginally) stable.

 This type of analysis has been started in \cite{PZ-CDC:14}, where partial results have been proposed showing that isotropic networks are difficult to be controlled, while anisotropic networks are more controllable. However, only intuitive definitions of isotropic and anisotropic networks were given there, supported by some illustrative examples. Actually, isotropic networks were described as networks in which there are no global preferential direction resulting from the network weights. Elaborating this idea, in the present paper we show that the right and left eigenvectors associated with the largest eigenvalue of the network adjacency matrix play an important role to determine whether the network is difficult to control. The importance of these eigenvectors is already well-known in network science, in connection with the concept of centrality, see \cite[Section 7.2]{Newman:10}. Precisely, we show that when these eigenvectors suggest that all the nodes have similar centrality degrees, then the network will be difficult to control. 

{\em \bfseries Paper structure.} The present paper is organized as follows. Section \ref{section:preliminary} is devoted to present some preliminary results. In Section \ref{sec:main-thm} we state the main result of the paper. Section \ref{sec:examples} contains some examples of application of the main theorem. Finally, in Section \ref{sec:concl} we summarize the contributions of the paper and we outline a number of possible future research directions. Most of the proofs of the more technical results are postponed to the appendix.

{\em \bfseries Notation.} We denote by $\R^{n}$ the set of $n$ dimensional vectors with real entries. The symbols $e_k$ and $\ones_n$ stand for the $k$-th vector of the canonical basis and the $n$-dimensional vector with all ones, respectively.
Moreover, we denote by $\R^{n\times m}$ the set of $n\times m$ matrices with real entries. The symbol $I_n$ stands for the $n\times n$ identity matrix. Given a matrix $A$, the symbol $A_{ij}$ means the $(i,j)$-th entry of $A$, while $A\T$ and $\ker(A)$ mean the transpose and the kernel of $A$, respectively. Finally, the symbol $\|A\|_{2}$ denotes the $2$-norm of $A$. 

A matrix $A\in\R^{n\times n}$ is said to be (Schur) stable if $\|A^t\|_{2}\to 0$ as $t\to\infty$, while it is said to be (Schur) marginally stable if $\|A^t\|_{2}$ is bounded in $t$. We denote the $i$-th component of the column vector $v\in\R^n$ as $v_i$. Given two vectors $v,w\in\R^{n}$, with $v\perp w$, we indicate that $v$ is orthogonal to $w$ in the standard Euclidean metric, that is $v\T w=0$. We let  $\diag(v)$, $v\in\R^n$, denote the $n\times n$ diagonal matrix with elements $v_1,\dots,v_n$ on the diagonal. Given a positive definite matrix $W\in\R^{n\times n}$,  we indicate with $\|x\|_{W}:=\sqrt{x\T W x}\,,\, x\in\R^n$, the weighted Euclidean norm (in case $W=I_{n}$ we simply write $\|x\|$). 

A matrix $A$ is said to be non-negative (positive) if $A_{ij}\geq 0$ ($A_{ij}> 0$) for all $i,j=1,\dots,n$. Irreducible and primitive matrices are special subclasses of non-negative matrices. Specifically, a  non-negative matrix $A\in\R^{n\times n}$ is said to be irreducible if $(I_n+A)^{k}$ is positive for some $k\in\N$, and primitive if $A^{k}$ is positive for some $k\in\N$ \cite[Chapter 8]{Horn-Johnson:85}.

Also, we denote by $\G=(\V,\E, A)$ the weighted graph with vertex (or node) set $\V=\left\{1,2,\dots,n\right\}$, edge set $\E\subseteq \V\times\V$ and adjacency matrix $A\in\R^{n\times n}$ satisfying $A_{ij}> 0$ if $(j,i)\in\E$. 
For a weighted graph $\G=(\V,\E, A)$ the out-degree of the node $i$ is equal to $\sum_iA_{ij}$, while the in-degree of the node $i$ is equal to $\sum_jA_{ij}$. A graph where the adjacency matrix has only {0-1} elements will be called unweighted. If $A=A\T$, the graph is called undirected. 
In case of undirected graphs, the in-degrees are equal to the out-degrees and they are called simply degrees and denoted by the symbol $\deg_A(i)$. In case the graph is undirected and unweighted, $\deg_A(i)$ corresponds to the number of edges of the node $i$.  

Finally, given two functions $f(x)$ and $g(x)$, with $g(x)$ non-zero, we write $f(x)=o(g(x))$ to mean that
\[
	\lim_{x\to\infty}\frac{f(x)}{g(x)}=0.
\]
Other standard notation is taken from \cite{PZB:14}.
\section{Preliminary facts}\label{section:preliminary}
Given a non-negative matrix $A\in\R^{n\times n}$, the following assumption will be used throughout the paper.

\begin{assum}\label{assum1}
Matrix $A$ is irreducible and marginally stable.
\end{assum}

These conditions on $A$ assure that there exists a unique real eigenvalue {\Rw  $\lambda_1$ such that $1\geq \lambda_1 \geq |\lambda_i|$ for all $i=2,\dots,n$,} where $\lambda_i$ denotes an eigenvalue of $A$; moreover, there exist positive left and right eigenvectors related to $\lambda_1$ \cite[Theorem 8.4.4]{Horn-Johnson:85}. 

One class of matrices which satisfies Assumption \ref{assum1} are irreducible column-stochastic (row-stochastic) matrices, i.e. irreducible matrices whose columns (rows) sum up to 1. This is an important class which has a large number of applications, e.g. in the context of distributed optimization and consensus algorithms \cite{Olfati:07, Garin-Schenato:11, Lynch:96}.

We denote with $v$ and $w$ the right and left eigenvector of $A$ with respect to $\lambda_1$ respectively, that is
\begin{equation}\label{eq:eigvec}
	Av=\lambda_1 v,\quad	A\T w=\lambda_1 w.
\end{equation}
The vectors $v$ and $w$ are called the right and left leading eigenvectors of the matrix $A$.
Since $v$ and $w$ are positive, we can define the following positive vector
\begin{equation}\label{eq:pi}
	\displaystyle	\pi	:=	\begin{bmatrix}\displaystyle \frac{w_1}{v_1} &\displaystyle \frac{w_2}{v_2}& \displaystyle \cdots & \displaystyle \frac{w_n}{v_n}\end{bmatrix}\T,
\end{equation}
and also the diagonal matrix 
\begin{align}\label{eq:Pi}
	\Pi:=\diag(\pi).
\end{align}
The vectors $v$ and $w$ can be normalized in such a way that $\sum_i\pi_i=1$.

We exploit matrix $\Pi$ to define an analogous to the time reversal of an irreducible Markov chain, see e.g. \cite[Section 1.6]{Levin:06}, that is we introduce the following matrix
\begin{align}\label{eq:AR}
	A^R:=\Pi\inv A\T \Pi.
\end{align}
This matrix is irreducible and $v$ and $w$ are respectively its right and left eigenvectors with respect to $\lambda_1$.
%
%
We can now define the matrix
\[
	A^S	:=	\Pi^{1/2}AA^R \Pi^{-1/2},
\]
which is symmetric and  positive semidefinite (and therefore has all real non-negative eigenvalues). Notice that, even if $A$ and $A^R$ are irreducible, $AA^R$ (and therefore also $A^S$) may not be irreducible. Let
\begin{equation}\label{eq:sigma2}
\sigma_{1} \geq \sigma_{2} \geq \dots \geq \sigma_{n} \geq 0,
\end{equation}
be the eigenvalues of $A^S$.
The following Lemma, whose proof is postponed to Appendix \ref{app:PrelFacts}, gives a characterization of the eigenvalues of $A^S$.

\begin{lem}\label{lm:eigvalAS}
The largest eigenvalue of $A^S$ is equal to the square of the largest eigenvalue of $A$, namely $\sigma_1=\lambda_1^{2}$.\footnote{Indeed this result holds for any irreducible $A$, independently from its stability.}
\end{lem}

The next Proposition, whose proof is given in Appendix \ref{app:PrelFacts}, is instrumental for Section \ref{sec:main-thm}.
\begin{prop}\label{cor:potenzeApery}
If $y\in\R^n$ is such that $y\T v=0$, then for any $t\geq 0$ {\Rw it holds}
\[
	\|(A\T)^t y\|_{\Pi\inv}^2 \leq \sigma_2^t \|y\|_{\Pi\inv}^2.
\]

\end{prop}

\section{The main result}\label{sec:main-thm}

Suppose that we are given a network represented by the graph $\G=(\V,\E,A)$, with  $A\in\R^{n\times n}$ satisfying Assumption \ref{assum1}\footnote{ We recall that irreducibility of $A$ implies that the graph $\G$ is connected.}. We consider the following discrete state-space system built according to the adjacency matrix of the graph:
\[
	x(t+1)=Ax(t)+B_\K u_{\K}(t),
\]
where $x\in\R^n$, $B_\K:=[e_{k_1}\, \cdots\, e_{k_m}]\in\R^{n\times m}$ and $u_{\K}\in\R^{m}$, with $\K:=\{k_1,\dots,k_m\}\subseteq \V$ being the set of control nodes and $m:=|\K|$.

{\Rw Our aim is to analyze the controllability of the system. However, we are not totally satisfied with the standard notion of controllability (from the origin) \cite{Kailath:80}:
\begin{defn}[Controllability] A system is controllable if, for every state $x_f\in\R^n$, there exists an input sequence $\left\{u(t)\right\}$ that steers the initial state $x(0)=0$ to $x_f$.
\end{defn}

In fact, we would also like to determine whether a controllable system  is difficult to control, where our idea of difficulty concerns the amount of energy needed to steer the initial state to a desired final one. To achieve this goal we introduce  the controllability Gramian of the system, which is defined as
}
\[
	\gram{T}:=\sum_{\tau=0}^{T-1} A^\tau B_\K B_\K\T (A\T)^\tau .
\] 

It is well known \cite{Kailath:80} that the system is controllable if and only if $\gram{T}$ is invertible for a big enough $T$. For controllable systems, denote by $u_{\K}^{*}(\cdot)$ the unique minimum $L^2$-energy input which steers the state from the initial value $x(0)=0$ to the final value $x(T)=x_f$. The energy of $u_{\K}^{*}(\cdot)$ is given by (see \cite{Kailath:80})
\[
	\mathrm{E}(x_f,T):=\sum_{\tau=0}^{T-1}\|u_\K^{*}(\tau)\|^2=x_f\T\gram{T}^{-1} x_f .
\]
From this observation it is possible to select a controllability metric. Indeed, choosing a worst-case analysis, we observe that the final state requiring the maximum energy is the one parallel to the eigenvector of $\gram{T}$ corresponding to its minimum eigenvalue, which is denoted by $\lmin{T}$. Precisely,
$$\max_{\|x_f\|=1}\mathrm{E}(x_f,T)=\lmin{T}^{-1}.$$
In this way, we can say that a controllable system is difficult to control if $\lmin{T}^{-1}$ is big, or equivalently if $\lmin{T}$ is small.

Finding upper bounds for $\lmin{T}$ can be therefore an instrument for proving that a system, although controllable in theory, is not controllable in practice due to the energy requirements.
The following Theorem provides a useful bound.

\begin{thm}
\label{mainThm}
Let $A$ be a non-negative matrix which satisfies Assumption \ref{assum1} and let the non-negative constant $\sigma_2$ and the positive {\Rw vectors $v$ and $w$ be defined as in (\ref{eq:sigma2}) and in (\ref{eq:eigvec}),} respectively.
If the product $AA\T$ is primitive, then
\begin{equation}\label{eq:eqthm}
	\lmin{T}\leq \frac{\max_i\{w_i/v_i\}}{\min_i\{w_i/v_i\}} \frac{\sigma_2^{\frac{n}{m}}}{\sigma_2^2(1-\sigma_2)}
\end{equation}
for all $T\in\N$ and all sets $\K$ with cardinality $m$.
\end{thm}

\begin{IEEEproof}
Given $\bar T\leq T$, the following chain of inequalities holds:
\begin{align*}
	& \lmin{\bar T} = \min_{\|x\|=1}x\T \gram{\bar T}x\\
	& =	\min_{\|x\|=1} x\T\left(\gram{T}-\sum_{\tau=\bar T}^{T-1}A^\tau B_\K B_\K\T (A\T)^\tau\right)x \\
	& \geq \lmin{T}+\min_{\|x\|=1} x\T \left(-\sum_{\tau=\bar T}^{T-1}A^\tau B_\K B_\K\T (A\T)^\tau\right) x \\
	& =   \lmin{T}-\max_{\|x\|=1}x\T \left(\sum_{\tau=\bar T}^{T-1}A^\tau B_\K B_\K\T (A\T)^\tau\right) x.
\end{align*}
We obtain in this way
\begin{align*}
	\lmin{T} 	& \leq \lmin{\bar T}\\
			& +\max_{\|x\|=1}x\T \left(\sum_{\tau=\bar T}^{T-1}A^\tau B_\K B_\K\T (A\T)^\tau\right) x.
\end{align*}
Let us now define
\[
	\lminhat{T}:=\min_{\substack{\|x\|=1\\x\perp v}} x\T \gram{T} x,
\]
where $v$ is as in equation (\ref{eq:eigvec}).
As before, if $\bar T\leq T$, we can prove that
\begin{align*}
	\lminhat{T} 	& \leq \lminhat{\bar T}\\
				& + \max_ {\substack{\|x\|=1\\x\perp v}} x\T\left(\sum_{\tau=\bar T}^{T} A^\tau B_\K B_\K\T (A\T)^\tau\right)x,
\end{align*}
and moreover it holds that $\lmin{T} \leq \lminhat{T}$.
Taking $\bar T=\left \lfloor{\frac{n-2}{m}}\right\rfloor$ we have that $\lminhat{\bar T}=0$. Indeed, $\gram{\bar T}$ coincides with $\Co{\bar T,\K}\Co{\bar T,\K}\T$, where $\Co{\bar T,\K}$ is the controllability matrix of the system, i.e.
\[
	\Co{\bar T,\K}:=\begin{bmatrix}
	B_\K & AB_\K & A^2B_\K & \cdots & A^{\bar T-1}B_\K
	\end{bmatrix}.
\]
Now, if we introduce the matrix $\Cob{\bar T,\K}:=\begin{bmatrix} \Co{\bar T,\K} & v \end{bmatrix}$ then the rank of $\Cob{t,\K}\Cob{t,\K}\T$ is strictly less than $n$. This in turn implies that there exists a $x\in\R^n$ such that $x\T \Cob{t,\K}\Cob{t,\K}\T x=0$ and hence $x$  is orthogonal to $v$ and is such that $x\T \gram{t}x=0$. 
From the fact that $\lminhat{\bar T}=0$, we have
\begin{align}
	& \lmin{T} 	\leq \lminhat{T} \notag \\
	& \quad\quad \leq \max_ {\substack{\|x\|=1\\x\perp v}} x\T\left(\sum_{\tau=\bar T}^{T-1}A^\tau B_\K B_\K\T (A\T)^\tau\right)x
					\le \beta_{\bar T}  \label{eq:beta}
\end{align}
where
\begin{align*}
\beta_{\bar T}:=\lim_{T\rightarrow\infty}\max_ {\substack{\|x\|=1\\x\perp v}} x\T\left(\sum_{\tau=\bar T}^{T-1}A^\tau B_\K B_\K\T (A\T)^\tau\right)x. 
\end{align*}
Our aim is to find a bound for $\beta_{\bar T}$. First observe that if we define $\alpha:=\min_i\left\{ v_i/w_i \right\}$, due to the particular form of $B_\K$, it holds that
\[
	y\T B_\K B_\K\T y\leq y\T y\leq \frac{1}{\alpha}y\T \Pi\inv y=\frac{1}{\alpha}\|y\|_{\Pi\inv}^2.
\]
Using Proposition \ref{cor:potenzeApery} and letting $y=(A\T)^\tau x$, we obtain
\[
	x\T A^{\tau}B_\K B_\K\T (A\T)^\tau x \leq \frac{1}{\alpha}\|(A\T)^\tau x\|_{\Pi\inv}^2
	\leq \frac{1}{\alpha}\sigma_2^{\tau}\| x\|_{\Pi\inv}^2.
\]
In this way we obtain
\[
	\beta_{\bar T}\leq \frac{1}{\alpha}  \max_{\substack{\|x\|=1 \\ x\perp v}} \|x\|_{\Pi\inv}^2\lim_{T\rightarrow\infty} \sum_{\tau=\bar T}^{T-1} \sigma_2^\tau.
\]
Now, since $\|x\|=1$, we have that
\begin{align*}
	\|x\|_{\Pi\inv}^2	&=	\sum_{i=1}^n \frac{v_i}{w_i} x_i^2\\
							&\leq	 \max_{i=1,\dots,n}\left\{ \frac{v_i}{w_i} \right\} \sum_{i=1}^n x_i^2 = \max_{i=1,\dots,n}\left\{ \frac{v_i}{w_i} \right\}.
\end{align*}
{\Rw Since $A$ is marginally stable and $AA\T$ is primitive, also $A^S$ is primitive (since it has the same sparsity pattern of $AA\T$). Therefore, it holds that $\sigma_2$ is strictly less than $\sigma_1 = \lambda_1^2\leq 1$, and therefore is strictly less than 1.}
Summing up all previous considerations, we have 
\[
	\beta_{\bar T}\leq \frac{\max_i\{v_i/w_i\}}{\min_i\{v_i/w_i\}}  \lim_{T\rightarrow\infty} \sum_{\tau=\bar T}^{T-1} \sigma_2^\tau=\frac{\max_i\{w_i/v_i\}}{\min_i\{w_i/v_i\}} \frac{\sigma_2^{\bar T}}{1-\sigma_2}.
\]
Letting $\bar T=\left \lfloor{\frac{n-2}{m}}\right\rfloor$, inequality (\ref{eq:beta}) yields $\lmin{T}\leq \beta_{\bar T}$. 
Since $\left \lfloor{\frac{n-2}{m}}\right\rfloor \geq n / m-2$ and $\sigma_2<1$, the statement follows.
\end{IEEEproof}

\begin{rem}
Concerning the primitivity condition on $AA\T$, a simple sufficient condition ensuring this property is $A$ irreducible and with all its diagonal elements strictly greater than zero. As a matter of fact, if $A_{ij}>0$, then 
\[
	(AA\T)_{ij}=A_{i1}A_{j1}+\dots+\underbrace{A_{ij}A_{jj}}_{>0}+\dots+A_{in}A_{jn}>0.
\]
This entails that all the elements greater than 0 in $A$ are greater than 0 also in $AA\T$, therefore $AA\T$ is irreducible and with positive diagonal elements, and thus primitive. This non-zero diagonal hypothesis is usually met by the networks we consider, since it implies that every node has its own dynamics. 
\end{rem}

\begin{rem}  \label{rem:stableMatr}
{\Rw In case $A$ is stable (and irreducible), the previous theorem holds true even without requiring the primitivity of $AA\T$. As a matter of fact, stability implies that $\lambda_1<1$ and it automatically follows that $\sigma_2\leq \sigma_1<1$.} It is worth noting that this version of our result is stronger than what it is proven in \cite[Theorem 3.1]{PZB:14}, namely that if $\|A\|_2<1$, then the system is certainly difficult to control. As a matter of fact a stable matrix does not necessarily have 2-norm less than one. For instance,  for all $\eps$ such that $0<\eps< \frac{3-\sqrt{5}}{2}$, the matrix $A_{\eps}:=\left[\begin{smallmatrix}\eps & 1 \\ \eps & \eps\end{smallmatrix}\right]$ is stable, but $\|A_{\eps}\|_{2}>1$.
\end{rem}

Once the system matrix $A$ and the number of controlled nodes $m$ are fixed, a metric that describes the controllability degree of the system is 
\[
	\Lambda(A,m)	:=\max_{\K\,:\, |\K|=m\atop T\in\N} \lmin{T}.
\]
Comparing several systems, those with a smaller $\Lambda(A,m)$ are more difficult to control, that is they need more energy for their control. The previous theorem shows that
\begin{equation}\label{eq:eqthm2}
	\Lambda(A,m)\leq \frac{\max_i\{w_i/v_i\}}{\min_i\{w_i/v_i\}} \frac{\sigma_2^{\frac{n}{m}}}{\sigma_2^2(1-\sigma_2)}.
\end{equation}

Now, consider a sequence of systems with matrices $A_n$ of increasing dimension $n$ and a number of controlled nodes $m(n)$ depending on $n$. In case
\begin{equation}\label{eq:limLambda}
\Lambda(A_n,m(n)) \xrightarrow{n\to\infty}0,
\end{equation}
we have that, even if for all $n$ the systems are controllable, for large $n$ the energy required to control them is so high that the control is practically impossible{\footnote  {\Rw Note that the tightness of bound \eqref{eq:eqthm} is not relevant when we consider the asymptotic behaviour as $n$ tends to infinity.}}. When the right-hand side of (\ref{eq:eqthm2}) tends to 0 as $n$ goes to infinity, (\ref{eq:limLambda}) is verified and so the systems become practically uncontrollable when $n$ is big.
Therefore, if we want to have a chance that the systems are practically controllable, we need to fix $m(n)$ in a way that the right-hand side of (\ref{eq:eqthm}) does not go to $0$ as $n$ increases (even though, since we have only an upper bound, this fact alone does not guarantee that the practically controllability can be in fact achieved).

To avoid that the right-hand side of (\ref{eq:eqthm2}) tends to 0 as $n$ goes to infinity, we can choose a fast enough growing $m(n)$. The growth of $m(n)$ will depend on two characteristics of $A$:
\begin{enumerate}
\item On how $\sigma_2$ depends on $n$: in some cases, as we will see in the next section, $\sigma_2$ stays bounded  away from $1$, namely there exists a constant $\varepsilon>0$ independent of $n$ such that $1-\sigma_2 \geq \varepsilon$, for all $n$. In other cases, instead, $\sigma_2$ tends to $1$ as $n$ goes to infinity. In this case it is important to understand how fast this occurs. The quantity $1-\sigma_2$, called the spectral gap of $A$, is essential to study inequality (\ref{eq:eqthm2}) in this case.
\item On how the fraction $\frac{\max_i\{w_i/v_i\}}{\min_i\{w_i/v_i\}}=\frac{\max_i\{\pi_i\}}{\min_i\{\pi_i\}}$ depends on $n$: as we will see in the following, in many important cases this fraction can be evaluated and in some of them this quantity remains bounded in $n$.
\end{enumerate} 
The way these two sources of dependence on $n$ interact will determine the growth of $m(n)$.
For instance, if $\sigma_2$ is bounded away from 1 and $\frac{\max_i\{w_i/v_i\}}{\min_i\{w_i/v_i\}}$ remains bounded in $n$, then practical controllability can be achieved only if the number of control nodes $m(n)$ grows linearly in $n$, namely only if a fixed fraction of the nodes are controlled. For example, this holds true for symmetric matrices \cite{PZ-CDC:14}. In fact, in this case, $\frac{\max_i\{v_i/w_i\}}{\min_i\{v_i/w_i\}}=1$. The following Proposition shows that this property holds for a more general class of matrices.

\begin{prop}\label{lm:condPi1} Let $v$ and $w$ be the right and left eigenvector of $A$, then
\begin{equation}\label{eq:lemmaB}
	\frac{\max_i\{w_i/v_i\}}{\min_i\{w_i/v_i\}}=1\Longleftrightarrow v\in\ker(A-A\T).
\end{equation}
\end{prop}

\begin{IEEEproof}
We have the following chain of equivalences
\begin{align*}
\frac{\max_i\{w_i/v_i\}}{\min_i\{w_i/v_i\}}=1 	& \Leftrightarrow v=\alpha w,\ \alpha\in\R\\
									&  \Leftrightarrow 	A\T v=\lambda_1 v\\				
									& \Leftrightarrow Av=A\T v \\
									& \Leftrightarrow (A-A\T )v=0 \\
									& \Leftrightarrow v\in\ker(A-A\T ) .		
\end{align*}
\end{IEEEproof}
The class of matrices satisfying \eqref{eq:lemmaB} comprises, for instance, irreducible doubly stochastic matrices, i.e. matrices which are both column- and row-stochastic.

\section{{\Rw Stochastic matrices, }Eigenvector Centrality and Reversible Matrices} \label{sec:centrality}
{\Rw The result of Theorem \ref{mainThm} acquires a nice interpretation when matrix $A$ is column-stochastic. Applying the theorem to such matrices, the next Corollary immediately follows:
\begin{cor} \label{cor:stoch}
If $A$ is an irreducible column-stochastic matrix  such that $AA\T$ is primitive, then for all $T\in\N$ and all sets $\K$ with cardinality $m$ it holds that
\begin{equation}\label{eq:eqthmcolsto}
\lmin{T}\leq \frac{v_{\max}}{v_{\min}} \frac{\sigma_2^{\frac{n}{m}}}{\sigma_2^2(1-\sigma_2)}
\end{equation}
with $v_{\max}:=\max_{i} v_{i}$ and $v_{\min}:=\min_{i} v_{i}$.\end{cor}

As a consequence of the Corollary{\footnote { Note that, with appropriate modifications, the Corollary can be applied also to row-stochastic matrices.}}
, the right eigenvector $v$, which also represents the invariant probability of the  Markov chain associated with the stochastic matrix $A$, plays a role in the controllability of the system represented by $A$. Vector $v$ also measures the (right) eigenvector centrality of the network associated with $A$ \cite[Chapter 7.2]{Newman:10}. It holds that, the bigger $v_i$ is, the more relevant or central the node $i$ is in the network, and therefore the fraction $\frac{v_{\max}}{v_{\min}}$ can be interpreted as a measure of heterogeneity in the node centralities.

Exploiting this interpretation of $v$, we have that a network where all the nodes have similar centrality has a lower heterogeneity index $\frac{v_{\max}}{v_{\min}}$ and, according to Corollary \ref{cor:stoch}, it will be more difficult to control\footnote{\Rw Interestingly, there is a connection between having a small heterogeneity index $\frac{v_{\max}}{v_{\min}}$ and the concept of wisdom of a network, as defined in \cite{Golub:10}.}.
On the other hand, easy to control networks need to have high heterogeneity in the nodes centrality.}

{\Rw The evaluation of the heterogeneity index $\frac{v_{\max}}{v_{\min}}$ is particularly simple if $A$ is reversible. The concept of reversibility for stochastic matrices is given in \cite{Levin:06}.}

\begin{defn}[Reversible matrix]
An  irredu\-cible  stochastic matrix $A\in\R^{n\times n}$ satisfying $A=A^R$ is called reversible.
\end{defn}

\begin{rem}
It is worth noticing that
\begin{itemize}
\item If the stochastic matrix $A$ is reversible, then $A$ can be regarded as the transition matrix of an irreducible time-reversible Markov chain as defined in \cite[Section 1.6]{Levin:06}.
\item If $A$ is reversible then it can be symmetrized by a diagonal transformation. Indeed, from (\ref{eq:AR}) it follows that
\[
	 S:= \Pi^{1/2} A \Pi^{-1/2}
\]
is symmetric, where $\Pi$ is defined in (\ref{eq:Pi}). This in turn implies that all the eigenvalues of $A$ are real, and therefore, by virtue of the irreducibility of $A$, they can be ordered in decreasing order as follows 
\[
	 1=\lambda_1 >\lambda_2\geq\cdots\geq \lambda_n\geq -\lambda_1.
\]
\item For a reversible matrix $A$, we have that $A_S=S^2$ and so, besides having $\sigma_1=\lambda_1^2$ as stated in Lemma \ref{lm:eigvalAS}, we also have that
\begin{equation}\label{eq:maxl2ln}
\sigma_2=\max\{\lambda_2^2,\lambda_n^2\}.
\end{equation}
\end{itemize}
\end{rem}
%

{\Rw It is possible to associate with a given non-negative matrix $C$ a stochastic matrix which is reversible if $C$ is symmetric.

Precisely, given an irreducible matrix $C\in\R^{n\times n}$, we can obtain a column-stochastic irreducible matrix $A$ as
\begin{equation}
	A := C\diag(\ones_n\T C)^{-1}.
	\label{eq:Pconstr}
\end{equation}
If $C=C\T$ then matrix $A$ is reversible\footnote{Notice that, even if $C$ is symmetric, $A$ is not necessarily so.}. The latter matrix represents the transition matrix of the (weighted) random walk built on the network represented by the original matrix. }

{\Rw Given a generic irreducible matrix $C$ (possibly not symmetric), the right leading eigenvector of $A$, obtained from $C$ as in \eqref{eq:Pconstr}, gives the (right) eigenvector centrality of the network associated with $A$, but it is also related to the PageRank \cite[Chapter 7.4]{Newman:10}, which is an important centrality measure concerning in this case the network associated with $C$\footnote{\Rw Indeed for a generic matrix $C$ the PageRank is defined as the right leading eigenvector of the column-stochastic matrix $A=\alpha C\diag(\ones_{n}\T C)^{-1}+(1-\alpha)\ones_{n}\ones_{n}^{\top}$, $0\leq\alpha\leq 1$. The factor $\alpha$, introduced to guarantee connectivity of the network, can be set to $1$ since in our case $C$ is irreducible by assumption.}.



Assuming now that $C$ is also symmetric, interesting analytic calculations can be carried out. Building $A$ as in \eqref{eq:Pconstr}, its leading right eigenvector $v$ has the following entries }
\begin{align}\label{eq:inv-meas-wrw}
	v_i = \frac{C^{\mathrm{col}}_i}{C_{\mathrm{tot}}},\ i=1,\dots,n,
\end{align}
with $C^{\mathrm{col}}_i:=\sum_j C_{ji}$  and $C_{\mathrm{tot}}:=\sum_{i,j} C_{ij}$. Consequently, the heterogeneity index is simply
\begin{equation}\label{eq:heterogeneity-index}
	\frac{v_{\max}}{v_{\min}}=\frac{\max_i\{C^{\mathrm{col}}_i\}}{\min_i\{C^{\mathrm{col}}_i\}}.
\end{equation}
Since $C^{\mathrm{col}}_i$ is the degree of the node $i$ in the undirected graph $\G=(\V,\E,C)$,
 Formula (\ref{eq:inv-meas-wrw}) says that in this case
the previously defined centrality coincides, up to a normalization, with the so-called ``degree centrality'' of the network represented by $C$  \cite[Chapter 7.1]{Newman:10}.

Assume now that $C$ is irreducible and symmetric and let
\begin{align*}
\underline{c}:=&\min\{C_{ij}\ |\ C_{ij}\not=0\} ,\\
\overline{c}:=&\max\{C_{ij}\ |\ C_{ij}\not=0\} .
\end{align*}
Moreover, assume that all the diagonal elements of $C$ are zero, i.e., $C_{ii}=0$, for all $i$.
In this case, it turns out that
\begin{equation}\label{eq:condPrw}
	\frac{v_{\max}}{v_{\min}}\leq \frac{(n-1)\overline{c}}{\underline{c}}.
\end{equation}
This shows that for reversible matrices obtained starting from matrices $C$ for which $\overline{c}/\underline{c}$ stays bounded, the fraction $v_{\max}/v_{\min}$ can grow at most linearly in $n$.

\section{ Examples} \label{sec:examples}

In this section, we show how the previous reasonings can be applied to three examples.
The first two examples concern random graph models, while the third one involves a structured graph.
Since these examples will deal with matrices with zero diagonal, in order to avoid unnecessary technical complications\footnote{In particular to avoid the presence of an eigenvalue in $-1$ and to also ensures that the product $AA\T$ is primitive.}, as in \cite{Levin:06} we use the concept of lazy version of a stochastic matrix. Precisely, from a column-stochastic matrix $A$, we can define a class of column-stochastic matrices as follows
\begin{equation}\label{eq:Palpha}
	A_{\alpha} := (1-\alpha) A+ \alpha I_{n}\, ,
\end{equation}
where $ 0<\alpha<1$. 
Notice that for any $\alpha$,
$A_{\alpha}$ has the same leading right eigenvector of $A$ and it is always primitive in case $A$ is irreducible.

\begin{exmp}[Weighted random walk on a B-A graph]\ \label{exmp:B-A}
The Barab\'asi-Albert preferential attachment graph (B-A for short), denoted in what follows by $\mathrm{BA}(n,d)$, with $d\geq 2$, is a well-known and widely studied random  graph model \cite[Chapter 14]{Newman:10}. Any realization of $\mathrm{BA}(n,d)$  is an undirected and unweighted graph with $n$ nodes and, in its construction, each newly added node makes $d$ connections with the previously existing nodes, according to a given rule that describes to which nodes the new node will be connected to (see \cite[Section 4.1]{Durrett:07}). An example of a graph constructed using this model can be found in Figure \ref{fig:BA-graph}.

\begin{figure}[h!]
\begin{center}
\includegraphics{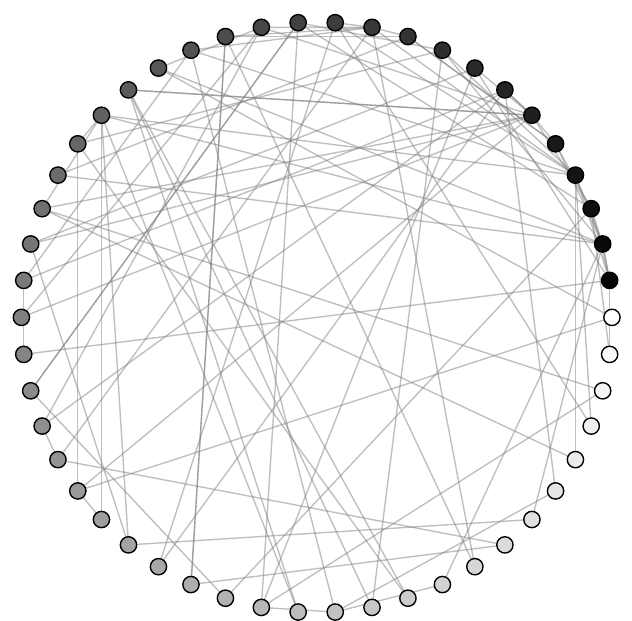}
\caption{ realization of a B-A graph $\mathrm{BA}(50,2)$ (the darker the node, the earlier it was inserted during the graph construction).}
\label{fig:BA-graph}
\end{center}
\end{figure}

In the following, we study the controllability degree of a lazy random walk on a weighted version of such a model. Specifically, given an adjacency matrix $\hat C$ of a B-A graph constructed as above, we associate with it a symmetric matrix $C$ obtained from $\hat C$ by letting $C_{ij}=0$ if $\hat C_{ij}=0$ and otherwise letting $C_{ij}$ and $C_{ji}$ equal to a number drawn uniformly and independently at random in the interval $[a, b]$, $0<a<b<\infty$.
Notice that the matrix $C$ defined in such a way is irreducible since the graph is connected. 
Using $C$, we can obtain matrix $A$ as in (\ref{eq:Pconstr}), and fixing a constant $\alpha \in ]0,1[$ we can build matrix $A_{\alpha}$ as in (\ref{eq:Palpha}). The latter represents a lazy weighted random walk on a B-A graph. In order to obtain information on the controllability degree of the network using Theorem \ref{mainThm}, we need to study the asymptotic behavior of $\sigma_2$.  Using Cheeger's inequality, it can be shown that the second eigenvalue of $A$, denoted with $\lambda_{2}(A)$, satisfies 
$$\lambda_2(A)\leq B$$
with high probability (w.h.p.) as  $n\rightarrow \infty$ (see Appendix \ref{app:weightedh} for the details), where $B\in]0,1[$ is a constant depending only on $a$ and $b$.
On the other hand, the $n$-th eigenvalue $\lambda_{n}(A)$ of $A$ satisfies 
$$\lambda_n(A)\geq -1.$$
These two facts imply, by (\ref{eq:maxl2ln}), that\footnote{Note that if $\lambda$ is an eigenvalue of $A$, $(1-\alpha)\lambda+\alpha$ is an eigenvalue of $A_\alpha$.}
\begin{align}
\sigma_2&=\max\{\lambda_2^{2}(A_\alpha),\lambda_n^2(A_\alpha)\}\notag\\
&\le\max_{z\in[-1,B]}\{((1-\alpha)z+\alpha)^2\}\notag\\
&=\max\{((1-\alpha)B+\alpha)^2,(2\alpha-1)^2\}<1.\label{eq:bound-lambda}
\end{align}
Thus, w.h.p. as $n\rightarrow \infty$, the number $\sigma_2$ is bounded by a constant smaller than one, depending only on $\alpha,a$ and $b$.
Now, since by (\ref{eq:condPrw}) $v_{\max}/v_{\min}\leq (n-1)b/a$, Corollary \ref{cor:stoch} yields
\[
	\Lambda(A_{\alpha},m)\leq (n-1)\frac{b}{a}\frac{\sigma_{2}^{\frac{n}{m}}}{\sigma_{2}^2(1-\sigma_{2})}.
\]
By taking logarithm of both sides in the latter expression and making further computations we can argue that
\begin{align}
	\log \Lambda(A_{\alpha},m)	& \leq k\log n  - \frac{n}{m}\log(1/ \sigma_{2}),
\end{align}
with $k>0$ being a constant depending only on $a$, $b$ and $\sigma_{2}$. Hence, since $\sigma_{2}$ is bounded away from $1$ w.h.p. as $n\rightarrow \infty$, the previous inequality allows us to conclude that, if $m=o(n/\log n)$, then
\[
	\log \Lambda(A_{\alpha},m)\xrightarrow{n\to\infty}-\infty \quad \text{ w.h.p.}
\]
or equivalently 
\[
	\Lambda(A_{\alpha},m)\xrightarrow{n\to\infty} 0 \quad \text{ w.h.p.}
\]
This implies that the weighted lazy random walks on almost all the realizations of B-A graph model are difficult to control by means of 
$o(n/\log n)$ driver nodes, as the cardinality of the network tends to infinity. 
\end{exmp}

\begin{exmp}[Weighted random walk on a E-R graph]\ \label{exmp:ER}
The Erd\"os-R\'enyi (E-R for short) graph model, denoted by $\mathrm{ER}(n,p)$, is one of the most celebrated and fundamental random graph models \cite[Chapter 11]{Newman:10}. A realization of the E-R graph model is an undirected and unweighted graph constructed as follows: starting from a graph of $n$ nodes, we place an undirected and unweighted edge between each distinct pair of nodes independently and with equal probability $p$. It can be shown that \cite[Theorem 2.8.1]{Durrett:07}, if $p=c\log n/n$, $c>1$, the realizations of a E-R graph model are connected w.h.p. as $n\rightarrow\infty$. With a partial abuse of language, we will refer to this subclass of E-R graph model as connected E-R graphs. Figure \ref{fig:ER-graph} shows an example of a connected E-R graph.

\begin{figure}[h!]
\begin{center}
\includegraphics{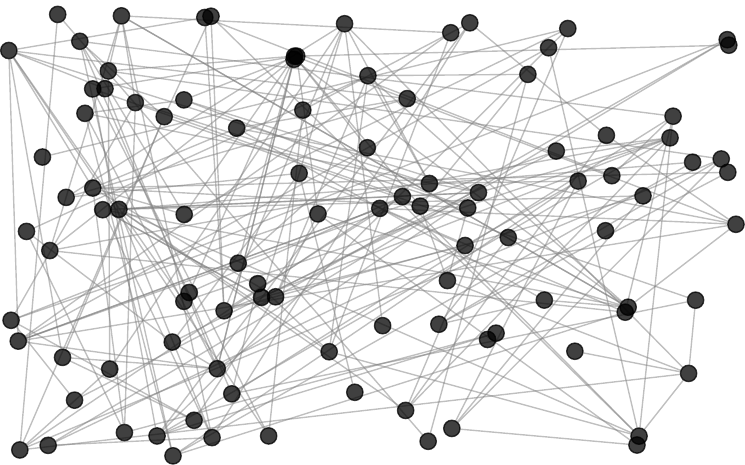}
\caption{A $100$ nodes connected E-R graph.}
\label{fig:ER-graph}
\end{center}
\end{figure}

Using the same procedure described in the previous example, we build a symmetric matrix $C$ with non-zero entries in $[a,b]$ from the adjacency matrix $\hat C$ of a connected E-R graph. From the analysis given in \cite[Section 6.5]{Durrett:07} and following the same lines of Example \ref{exmp:B-A}, we can argue that for the weighted lazy random walk on the graph described by $C$, the number $\sigma_{2}$ is bounded away from $1$ w.h.p. as $n$ tends to infinity. Hence, we are in position to conclude that, by virtue of Corollary \ref{cor:stoch}, the weighted lazy random walks on almost all connected E-R graphs are difficult to control w.h.p. as $n$ goes to infinity if we use $m=o(n/\log n)$ control nodes.
\end{exmp}

\begin{rem}
Interestingly, the choice made for the probability $p$ allows us to make, in this case, further considerations on the heterogeneity index $v_{\max} / v_{\min}$. As a matter of fact, for an unweighted E-R graph with $p=c\log{n}/n$, it holds that the maximum and minimum degree (denoted by $\deg_{\max}$ and $\deg_{\min}$, respectively)  satisfy \cite[Ex. 3.4]{Bollobas:01} 
\begin{align*}
\deg_{\max} \xrightarrow{n\to\infty} \overline{\gamma} \log n \quad \text{ w.h.p.}\\
\deg_{\min} \xrightarrow{n\to\infty}  \underline{\gamma} \log n \quad \text{ w.h.p.}
\end{align*}
where $0<\underline{\gamma}<\overline{\gamma}<\infty$ depend only on $c$. Consequently, since by (\ref{eq:heterogeneity-index}) we have that
\[
	\frac{v_{\max}}{v_{\min}}\leq\frac{{\deg_{\max}b}}{{\deg_{\min}a}},
\] 
we can argue that $v_{\max}/v_{\min}$ is w.h.p. upper bounded by a constant as $n\to\infty$.
From this fact we can conclude that a stronger result holds, namely that this class of systems
are difficult to control w.h.p. as $n$ goes to infinity, even if we use $m=o(n)$ control nodes.
\end{rem}

\begin{figure*}[!ht]
\begin{center}
\includegraphics[scale=1]{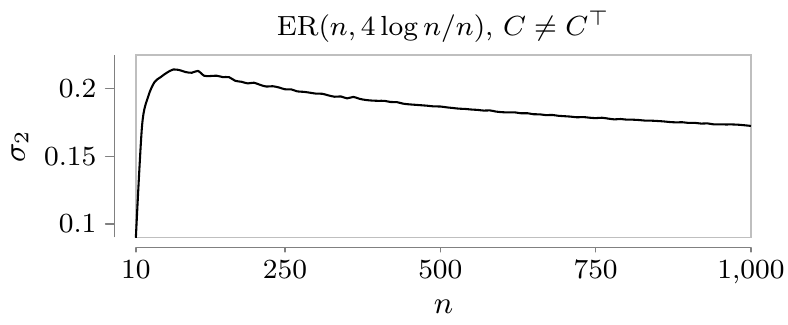}
\hspace{0.75cm}
\includegraphics[scale=1]{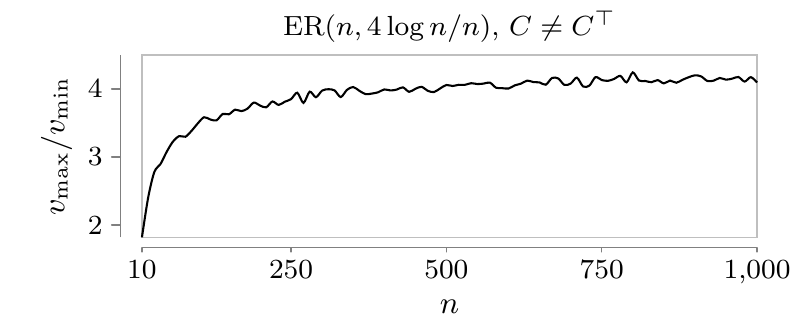}

\vspace{0.2cm}

\includegraphics[scale=1]{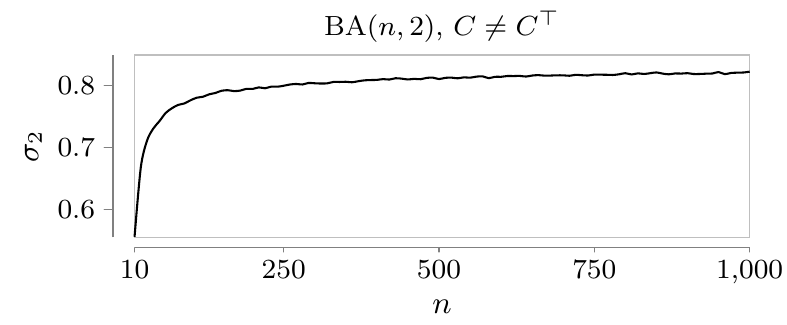}
\hspace{0.75cm}
\includegraphics[scale=1]{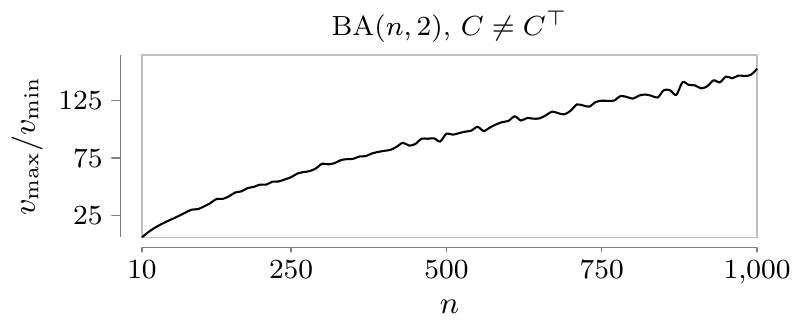}

\caption{The asymptotic behaviour of $\sigma_{2}$ and $v_{\max}/v_{\min}$ for networks with asymmetric matrix $C$, $\alpha=1/2$ and $a=1/2$, $b=4$. We start from an undirected random graph and, we associate with each edge $(i,j)$ of the graph a weight $C_{ij}$ drawn uniformly and independently at random in the interval $[a, b]$. The upper figures refer to E-R graphs $\mathrm{ER}(n,4\log n/n)$, while lower figures to B-A graph $\mathrm{BA}(n,2)$. For each $n$, values are obtained averaging 250 realizations of each random model. The simulations were performed using Python's NetworkX library \cite{hagberg:08}.} 
\label{fig:simulations-asymmetric}
\end{center}
\end{figure*}

\begin{rem}
The symmetry of the matrix $C$ in Formula (\ref{eq:Pconstr}) is necessary to enable a mathematical proof of the results proposed in the previous examples. In fact, if $C$ is not symmetric, we have no mathematical instruments  to estimate $\sigma_{2}$ and $v_{\max}/v_{\min}$. 
One could wonder how much this condition is crucial to observe practical uncontrollability of these networks. To investigate this issue, we carried out some simulations in which we did not impose symmetry on $C$. Figure \ref{fig:simulations-asymmetric} shows that for B-A random graphs $\sigma_{2}$ stays bounded away from $1$ and $v_{\max}/v_{\min}$ grows at most linearly in $n$, while  for E-R random graphs $\sigma_{2}$ stays bounded away from $1$ and $v_{\max}/v_{\min}$ stays bounded.
From this one can  infer that symmetry seems not to be a crucial property and that the systems built in this way remains difficult to be controlled.
\end{rem}

\begin{exmp}[Weighted random walk on a cube]\
Consider the unweighted graph consisting in the $3$-dimensional $k$-ary array \cite{Azizoglu:99}, defined as the Cartesian product of three path graphs of length $k$ (see Fig. \ref{fig:cube-graph}). In analogy with the previous examples, starting from the adjacency matrix $\hat C$ of this graph, we build a symmetric matrix $C$ having non-zero entries in $[a,b]$. Let $A\in\R^{n\times n}$, $n:=k^{3}$, be obtained from $C$ as in the previous examples and denote by $A_\alpha$ the lazy version of $A$.

\begin{figure}[h!]
\begin{center}
\includegraphics{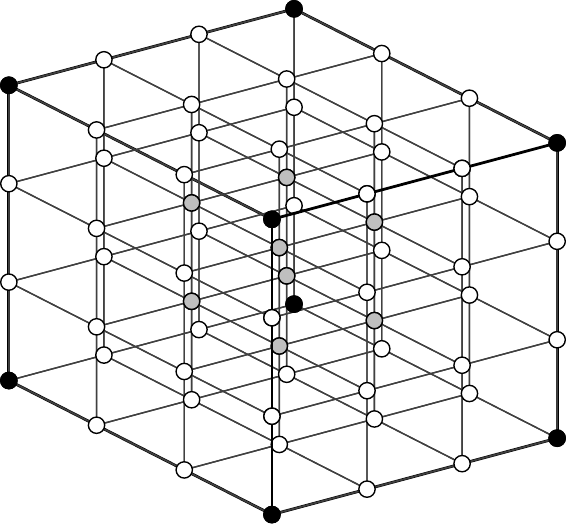}
\caption{$3$-dimensional $4$-ary array:  corner nodes are filled in black, while internal nodes in gray. 
}
\label{fig:cube-graph}
\end{center}
\end{figure}

By (\ref{eq:heterogeneity-index}), we have that
\[
	\frac{v_{\max}}{v_{\min}} \leq{ \frac{2b}{a }}
\]
since the maximum degree for the unweighted graph is $6$ (for internal nodes) and the minimum is $3$ (for corner nodes). 
 In order to estimate $\lambda_2(A)$ we exploit again Cheeger's inequality, which yields that (see Appendix \ref{app:weightedh})
\begin{equation}\label{eq:cbcube}
	\lambda_2(A)\leq 1 - K n^{-2/3}=:B(n),
\end{equation}
where $K$ is a positive constant depending only on $a$ and $b$. Now, using (\ref{eq:bound-lambda}), it follows that
\begin{align*}
	\sigma_2		& =		\max\{\lambda_2^{2}(A_\alpha),\lambda_n^2(A_\alpha)\}\\
				& \leq 	\max\{((1-\alpha)B(n)+\alpha)^2,(2\alpha-1)^2\}.
\end{align*}
Notice that, for $n$ sufficiently large, it holds that
\[
\max\{((1-\alpha)B(n)+\alpha)^2,(2\alpha-1)^2\}=	 ((1-\alpha)B(n)+\alpha)^2.
\]
The latter inequality can be used in (\ref{eq:eqthm2}), obtaining (see Appendix \ref{app:weightedh} for detailed calculations)
\begin{align*}
	\log(\Lambda(A_{\alpha},m)) 	& \leq 	k_0 	-k_{1} \frac{n^{1/3}}{m}+k_{2}\log n
\end{align*}
where $k_0$, $k_1$ and $k_{2}$ are positive constant values which depend only on $a$, $b$ and $\alpha$. 
When $n$ tends to infinity we have that, as long as we use a number $m=o(n^{1/3}/\log n)$ 
of control nodes, the system is practically uncontrollable since $\log(\Lambda(A_{\alpha},m))\rightarrow -\infty$, and therefore $\Lambda(A_{\alpha},m)\rightarrow 0$.
\label{exmp:cube}
\end{exmp}

{\Rw
\begin{rem}
In a similar way as Example \ref{exmp:cube}, it is possible to prove that, given $d\geq3$, the $d$-dimensional $k$-ary array is practically uncontrollable with a portion $m=o\left(\frac{n^{1-\frac{2}{d}}}{\log n}\right)$ of controllers. In case $d=1$ or $d=2$, our bound is not useful to determine the controllability degree of the system.
\end{rem}
}

\begin{figure*}[!ht]
\begin{center}
\includegraphics[scale=1]{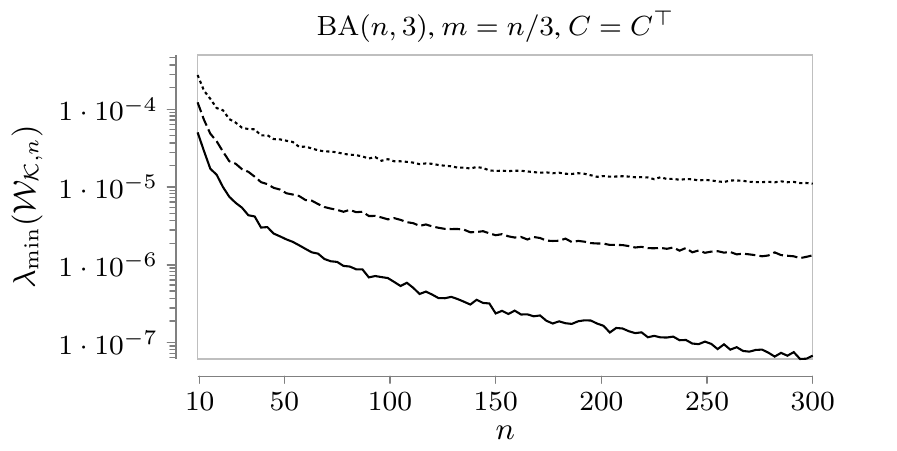}
\hspace{-0.35cm}
\includegraphics[scale=1]{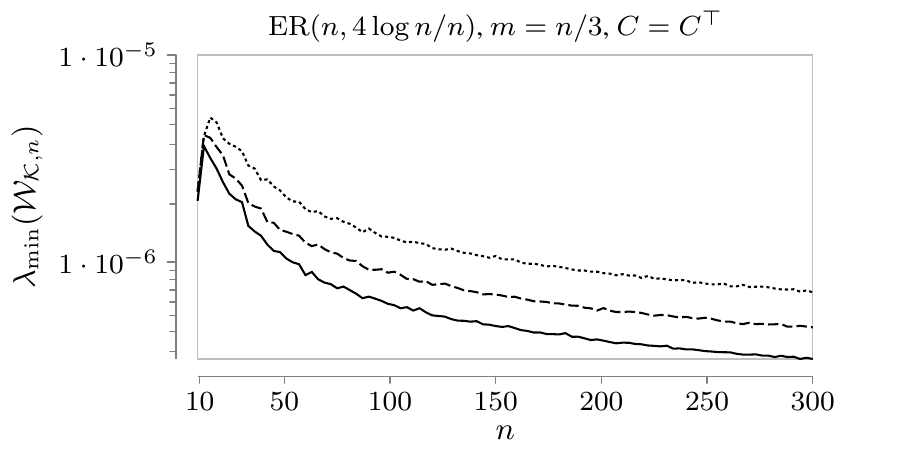}

\vspace{0.1cm}
\includegraphics[scale=1]{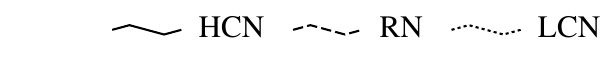}
\hspace{2.6cm}
\includegraphics[scale=1]{./legend3}

\caption{Plots of the controllability degree $\lmin{n}$ for B-A and E-R networks with symmetric matrix $C$, $\alpha=1/2$, $a=1/2$, $b=2$ and using $m=n/3$ controllers selected in three different ways. Precisely, the controllers are placed at the nodes with highest PageRank centrality (Highest Centrality Nodes = HCN), with lowest centrality (Lowest Centrality Nodes = LCN), and at nodes selected at random (Random Nodes = RN). For each $n$, values are obtained averaging 500 realizations of each random model.} 
\label{fig:simulations-positionning}
\end{center}
\end{figure*}

\section{Conclusions and future directions}\label{sec:concl}

In this paper we have shown the relevance of the network centrality for understanding how much energy is needed to control a dynamical network with adjacency matrix which is irreducible and (marginally) stable. Precisely, it is shown that if the matrix describing such a linear dynamical network has rather uniform node centralities, then it will be difficult to control. 

There are still many open questions that need to be addressed. One is related to whether or not also the converse of the result shown in this paper holds true, namely if a system characterized by nodes with very different centralities is easy to be controlled. Simulations seem to disprove this conjecture, since node centrality alone seems not to be enough for characterizing the network controllability.

Moreover, another interesting aspect concerns if node centralities could be helpful in selecting control nodes, in case we want to increase the controllability degree. Indeed, one could argue that it is more convenient to control the nodes with largest centralities. Figure \ref{fig:simulations-positionning} contains some results obtained carrying out some simulations in order to clarify this question. In these figures it is shown the value of $\lmin{T}$ as a function of $n$ in three cases: random positioning of the control nodes, control nodes positioned at the nodes with largest centrality and control nodes positioned at the nodes with smallest centrality. {\Rw Both for E-R and B-A random graphs something quite unexpected happen. Indeed in both cases positioning the control nodes at the nodes with smallest centrality happens to be the best strategy. We are unable to interpret this counterintuitive behavior{\footnote {\Rw Note that there seems to be a connection with the results found in \cite{Liu:11} concerning which are the most important nodes to control.}}, whose understanding surely deserves our attention in our future work on this subject.}


Another still open question is to prove what it seems quite intuitive, namely that if we can control a fixed percentage of the nodes, the energy needed to control the system stays bounded, or, more formally, if $m=a n$, where $a$ is any number in $]0,1[$, then $\Lambda(A,m)$ stays bounded away from zero.

These are only a few questions that need to be better understood on this subject which, although it has received a lot of attention in the last years, still remains very challenging and full of issues that require to be investigated.


\appendices

\section{Proofs and additional results}

\label{app:PrelFacts}

\begin{IEEEproof}[Proof of Lemma \ref{lm:eigvalAS}]
We first note that the positive vector
\[
	r	:= \Pi^{1/2}v= \begin{bmatrix} \sqrt{v_1 w_1} & \cdots & \sqrt{v_n w_n} \end{bmatrix}\T 
\]
 is such that $A^S r=\lambda_1^2 r$, as can be seen by direct computation. Since $A^S$ is a non-negative matrix and $r$ is positive, we can use \cite[Corollary 8.1.30]{Horn-Johnson:85} to conclude that $\lambda_1^2$ is the spectral radius of $A^S$. Therefore $\sigma_1=\lambda_1^2$.
 \end{IEEEproof}
 
\begin{IEEEproof}[Proof of Proposition \ref{cor:potenzeApery}] We first prove two instrumental lemmas.
\begin{lem}\label{lm:1}
If $x\in\R^n$ is such that $x\T w=0$, then 
\[
	\|A^R x\|_\Pi^2\leq \sigma_2 \|x\|_\Pi^2,
\]
where $\sigma_2$ is the second largest eigenvalue of $A^S$.
\end{lem}

\begin{IEEEproof}
First note that 
\[
	\|A^R x\|_\Pi^2 = x\T (A^R)\T \Pi A^R x = x\T \Pi^{1/2}A^S \Pi^{1/2}x. 
\]
Since $x\T w=0$, we also have that $x\T \Pi^{1/2}$ is orthogonal to $r$, which is the eigenvector of $A^S$ with respect to $\sigma_1$. Therefore
\[
	x\T \Pi^{1/2} A^S\Pi^{1/2}x \leq \sigma_2 x\T \Pi^{1/2}\Pi^{1/2}x=\sigma_2 \|x\|_\Pi^2
\]
and this ends the proof. 
\end{IEEEproof}

\begin{lem}\label{lm:normaPiInvy}
If $y\in\R^n$ is such that $y\T v=0$, then 
\[
	\|A\T y\|_{\Pi\inv}^2\leq \sigma_2 \|y\|_{\Pi\inv}^2.
\]
\end{lem}

\begin{IEEEproof}
Defining $x:=\Pi\inv y$, we have these two facts: 
\begin{align*}
	& x\T w  = y\T \Pi\inv w=y\T v=0\,,\\
	& \|x\|_\Pi^2= x\T\Pi x=y\T\Pi\inv\Pi\Pi\inv y=\|y\|_{\Pi\inv}^2.
\end{align*}
Moreover
\begin{align*}
	\| A^R x\|_\Pi^2	&=	x\T \Pi A \Pi\inv \Pi \Pi\inv A\T\Pi x \\
							&=	y\T A\Pi\inv A\T y \\
							&=	\| A\T y\|_{\Pi\inv}^2.
\end{align*}
Exploiting the previous lemma we have
\[
	\|A\T y\|_{\Pi\inv}^2	=	\| A^R x\|_\Pi^2 \leq \sigma_2 \|x\|_\Pi^2=\sigma_2 \|y\|_{\Pi\inv}^2.
\]
and we are done. 
\end{IEEEproof}

Using the lemmas we have just proven, we can finally prove the proposition as follows.
Since $y\T A^{t-1}v=\lambda_1^{t-1}y\T v=0$, defining $x:=(A\T)^{t-1}y$, we have that $x\T v=0$. Then
\begin{align*}
	\|(A\T)^ty\|_{\Pi\inv}^2	&=	y\T A^t \Pi\inv (A\T)^t y \\
									&=	x\T A \Pi\inv A\T x \\
									&=	\|A\T x\|_{\Pi\inv}^2 \\
									&\leq \sigma_2 x\T \Pi\inv x \\
									&=	\sigma_2\|(A\T)^{t-1}y\|_{\Pi\inv}^2
\end{align*}
and applying the same reasoning $t$ times, we obtain
\[
	\|(A\T)^t y\|_{\Pi\inv}^2\leq\sigma_2^t\|y\|_{\Pi\inv}^2\,. 
\]
\end{IEEEproof}

\section{Detailed Calculations for Examples \ref{exmp:B-A} and \ref{exmp:cube}}\label{app:weightedh}

First observe that, due to the Cheeger inequality \cite[Section 6.2]{Durrett:07}
\begin{align}
	\lambda_2(A)	&\le	 1-\frac{h^2}{2},
\label{eq:CheegerBound}
\end{align}
where $h>0$ is known as the bottleneck ratio or Cheeger constant of the column-stochastic matrix $A$. This constant is defined as
\begin{equation}\label{eq:conductance}
h=\min_{\substack{S\subset \V \\ v(S)\leq 1/2}}\frac{Q(S,\bar{S})}{v(S)} ,
\end{equation}
where $S$ is a subset of $\V$, $\bar{S}$ is its complement, $v(S):=\sum_{k\in S}v_k$, and 
$$Q(S,\bar{S}):=\sum_{i\in S,j\in \bar{S}} A_{ji}v_{i}.$$
Since in our case the column-stochastic matrix $A$ is defined from a symmetric matrix $C$ as in (\ref{eq:Pconstr}), then 
\begin{equation}\label{eq:estQ}
	Q(S,\bar{S})= \frac{1}{ C_{\textrm{tot}}}\sum_{i\in S,j\in \bar{S}}  C_{ji}\geq \frac{a}{C_{\textrm{tot}}} e(S,\bar{S})
\end{equation}
where $e(S,\bar{S})$ denotes the number of edges between $S$ and $\bar{S}$ in the unweighted graph described by $\hat{C}$.
The denominator of (\ref{eq:conductance}) can be treated as follows
\begin{equation}\label{eq:estvS}
v(S)=\sum_{k\in S}v_k=\sum_{k\in S} \frac{C^{\text{col}}_k}{C_\text{tot}}\leq \frac{b}{C_\text{tot}}\mathrm{vol}(S)
\end{equation}
where $\mathrm{vol}(S)$ is the so-called volume of $S$ in the unweighted graph and coincides with the total number of edges of the nodes belonging to $S$.
From the previous inequalities we can argue that
\[
h\geq \frac{a}{b}\min_{\substack{S\subset \V \\ v(S)\leq 1/2}}\frac{e(S,\bar{S})}{\mathrm{vol}(S)}.
\]

In \cite[Section 6.4]{Durrett:07} it is shown that, for the B-A random graph, the minimum in the previous formula is lower bounded by a positive constant w.h.p. as $n\to\infty$. This shows that also $h$ is lower bounded by a positive constant w.h.p. as $n\to\infty$. This solves Example \ref{exmp:B-A}.

As far as Example \ref{exmp:cube} is concerned, we still use the Cheeger bound (\ref{eq:CheegerBound}) and the inequalities (\ref{eq:estQ}) and (\ref{eq:estvS}). Then, since in this case $\mathrm{vol}(S)\le 6|S|$, we can argue that 
\[
	\frac{Q(S,\bar{S})}{v(S)}\geq \frac{a}{6b}\frac{e(S,\bar{S})}{|S|}
\]
which implies that
\begin{align*}
	h\geq	 \frac{a}{{6} b}\min_{\substack{S\subset \V \\|S|\leq n/2}} \frac{ e(S,\bar{S})}{|S|},
\end{align*}
where the minimum of $e(S,\bar{S})/|S|$ is the so-called isoperimetric number of the graph \cite[Chapter 3]{Chung:97}, which in a $3$-dimensional $k$-ary array is lower bounded by ${2}/{n^{1/3}}$ \cite{Azizoglu:99}. Therefore in this case we have
\[
	h\geq \frac{a}{\displaystyle { 3} b}\frac{1}{n^{1/3}}.
\] 
Now, by virtue of the Cheeger inequality (\ref{eq:CheegerBound}), it holds that
\[
	\lambda_{2}(A) \leq 1-Kn^{-2/3}=:B(n),
\]
where $K:=\frac{a^{2}}{18 b^{2}}$ is positive.

Concerning the evaluation of $\Lambda(A_{\alpha},m)$, by Corollary \ref{cor:stoch} and inequality (\ref{eq:eqthm2}), it follows that
\[
	\Lambda(A_{\alpha},m) \leq {\frac{2 b }{a}} \frac{\bar{B}(n)^{{\frac{n}{m}}}}{\bar{B}(n)^{2}\left(1-\bar{B}(n)\right)},
\]
where 
$$\bar{B}(n):=((1-\alpha)B(n)+\alpha)^2=(1-(1-\alpha)Kn^{-2/3})^2.$$
Now take logarithm of both sides 
\begin{align*}
	\log(\Lambda(A_{\alpha},m)) 	& \leq 	\log{\frac{2 b }{a}} + \frac{n}{m}\log\bar{B}(n) \\
	&-2 \log\bar{B}(n) -\log (1-\bar{B}(n)).\end{align*}
Using the Taylor's expansion as $n\to\infty$ we get
\begin{align*}
	\log\bar{B}(n) 	& =2\log(1-(1-\alpha)Kn^{-2/3})\\
							& \simeq	-2(1-\alpha)Kn^{-2/3}
\end{align*}
and
\begin{align*}
	\log(1-\bar{B}(n)) 	& =\log(2(1-\alpha)Kn^{-2/3}-(1-\alpha)^2K^2n^{-4/3})\\
							& \simeq	\log(2(1-\alpha)Kn^{-2/3})\\
							& =	\log(2(1-\alpha)K)-\frac{2}{3}\log n.
\end{align*}
Considering finally that $-2\log\bar{B}(n)\le 4(1-\alpha)K$ we can argue that
\begin{align*}
	\log(\Lambda(A_{\alpha},m)) 	& \leq 	k_0 	-k_{1} \frac{n^{1/3}}{m}+k_{2}\log n,
\end{align*}
where $k_0$, $k_1$ and $k_{2}$ are positive constant values which depend only on $a$, $b$ and $\alpha$.  
%
%

\bibliographystyle{IEEEtran}
\bibliography{IEEEabrv,Bibliography}

\end{document}